# A cardiac electromechanics model coupled with a lumped parameters model for closed-loop blood circulation. Part I: model derivation

Francesco Regazzoni[1,*], Matteo Salvador[1,*], Pasquale Claudio Africa[1], Marco Fedele[1], Luca Dedè[1], Alfio Quarteroni[1,2]

[1] MOX-Dipartimento di Matematica, Politecnico di Milano, Milan, Italy
[2] Professor Emeritus, École polytechnique fédérale de Lausanne, Lausanne, Switzerland

* These authors equally contributed to this work

**Abstract**

We propose an integrated electromechanical model of the human heart, with focus on the left ventricle, wherein biophysically detailed models describe the different physical phenomena concurring to the cardiac function. We model the subcellular generation of active force by means of an Artificial Neural Network, which is trained by a suitable Machine Learning algorithm from a collection of pre-computed numerical simulations of a biophysically detailed, yet computational demanding, high-fidelity model. To provide physiologically meaningful results, we couple the 3D electromechanical model with a closed-loop 0D (lumped parameters) model describing the blood circulation in the whole cardiovascular network. We prove that the 3D-0D coupling of the two models is compliant with the principle of energy conservation, which is achieved in virtue of energy-consistent boundary conditions that account for the interaction among cardiac chambers within the computational domain, pericardium and surrounding tissue. We thus derive an overall balance of mechanical energy for the 3D-0D model. This provides a quantitative insight into the energy utilization, dissipation and transfer among the different compartments of the cardiovascular network and during different stages of the heartbeat. In virtue of this new model and the energy balance, we propose a new validation tool of heart energy usage against relationships used in the daily clinical practice. Finally, we provide a mathematical formulation of an inverse problem aimed at recovering the reference configuration of one or multiple cardiac chambers, starting from the stressed configuration acquired from medical imaging. This is fundamental to correctly initialize electromechanical simulations. Numerical methods and simulations of the 3D-0D model will be detailed in Part II.

## 1 Introduction

In this two-part series of papers, we present novel mathematical and numerical models of cardiac electromechanics [1, 2, 7, 18, 28, 32, 34, 56]. Part I is devoted to present and analyze the mathematical models; Part II [45] deals with the numerical methods and simulations. Specifically, we present a fully coupled model of the cardiac function, which results from the concerted action of several physical phenomena – electrophysiology, biochemistry, mechanics, fluid dynamics – interacting at different spatial and temporal scales [10, 15, 54, 55]; these range from nanometers to centimeters and from nanoseconds to seconds, respectively [5, 26]. We describe all of these phenomena in terms of biophysically detailed models, written as systems of PDEs (partial differential equations) and ODEs (ordinary differential equations), which realize the coupling of elecrophysiology, ionic and gating models, active force generation at the cellular and tissue levels, passive mechanics of the muscle and blood flow in the chambers and in the circulatory system. Among the original contributions of this paper, we present a 0D (zero-dimensional) hemodynamics model of the whole cardiovascular system, which is coupled with our 3D electromechanical model of specific cardiac chambers to form a closed-loop 3D-0D circulation model. Such coupling is an essential step to provide physically meaningful numerical simulations, which are performed considering a 3D electromechanical description of the left ventricle only.

Among the different physical phenomena concurring to the heart function, the intrinsic complexity of the subcellular mechanisms underlying active force generation typically dictates the use of phenomenological models [29, 33, 46] or Monte Carlo approximations of physics-based models [23, 59, 60], which are however characterized by large computational costs. In our electromechanical model, instead, we rely on a physics-based model of force generation, with explicit representation of the cooperative interactions among the subcellular units [41]. As this model features more than 2000 internal variables, we surrogate it thanks to the Machine Learning algorithm that we proposed in [44]. Specifically, we train – in an offline phase – a reduced model, based on Artificial Neural Networks (ANNs), that approximates within a prescribed accuracy the results of the original one. By using this ANN-based model (featuring only two internal variables) instead of the high-fidelity one in the multiscale electromechanical simulation, we reduce the number of variables from more than 2000 to only 2, with an overall relative error of the order of $10^{-3}$ [44].

We prove that our closed-loop circulation model satisfies a balance of mechanical energy. Through the calculation of the different terms of this balance during an heart beat, we provide a quantitative insight in the cardiac energy distribution, highlighting the features of different compartments and the different stages of the heartbeat, i.e. when energy is injected, dissipated or transformed. Thanks to our model we can also assess the validity of simplified relationships commonly used in the clinical practice to estimate the main indicators of heart energy distribution [25].

We prove that the coupling between the 3D electromechanical model and the 0D circulation model is consistent with the principles of energy conservation. Indeed, we impose a boundary condition at the base of the left ventricle that we purposely denote as energy-consistent boundary condition [40, 44]. Moreover, we apply a boundary condition at the epicardium that keeps into account the interaction between the left ventricle and the pericardial sac [19, 54].

Looking towards the patient-specific customization of our model, we note that cardiac geometries are always acquired in vivo through imaging techniques, for which the left ventricle is loaded, mainly by the pressure acting on the endocardium. Therefore, the stress-free configuration, to which the equations for cardiac electromechanics must refer, is not known a priori. As this is necessary to set the reference configuration for the mechanical model, we formulate an inverse problem aimed at recovering such stress-free reference configuration starting from the geometry acquired from medical imaging, whose solution will be performed in Part II [45].

## 1.1 Paper outline

This paper is organized as follows. In Sec. 2 we describe a mathematical model for cardiac electromechanics endowed with our biophysically detailed ANN-based active stress model, and a lumped parameters model for the circulation system. Then, in Sec. 3 we address the problem of initializing the simulation on a physically sound basis, by formulating an inverse problem aimed at recovering the stress-free configuration from a stressed one. In Sec. 4 we derive a balance of the mechanical energy for both the 0D blood circulation model and the 3D-0D electromechanical-circulation coupled model and we provide a quantitative analysis of the different terms involved in the balance. Finally, in Sec. 5, we draw our conclusions.

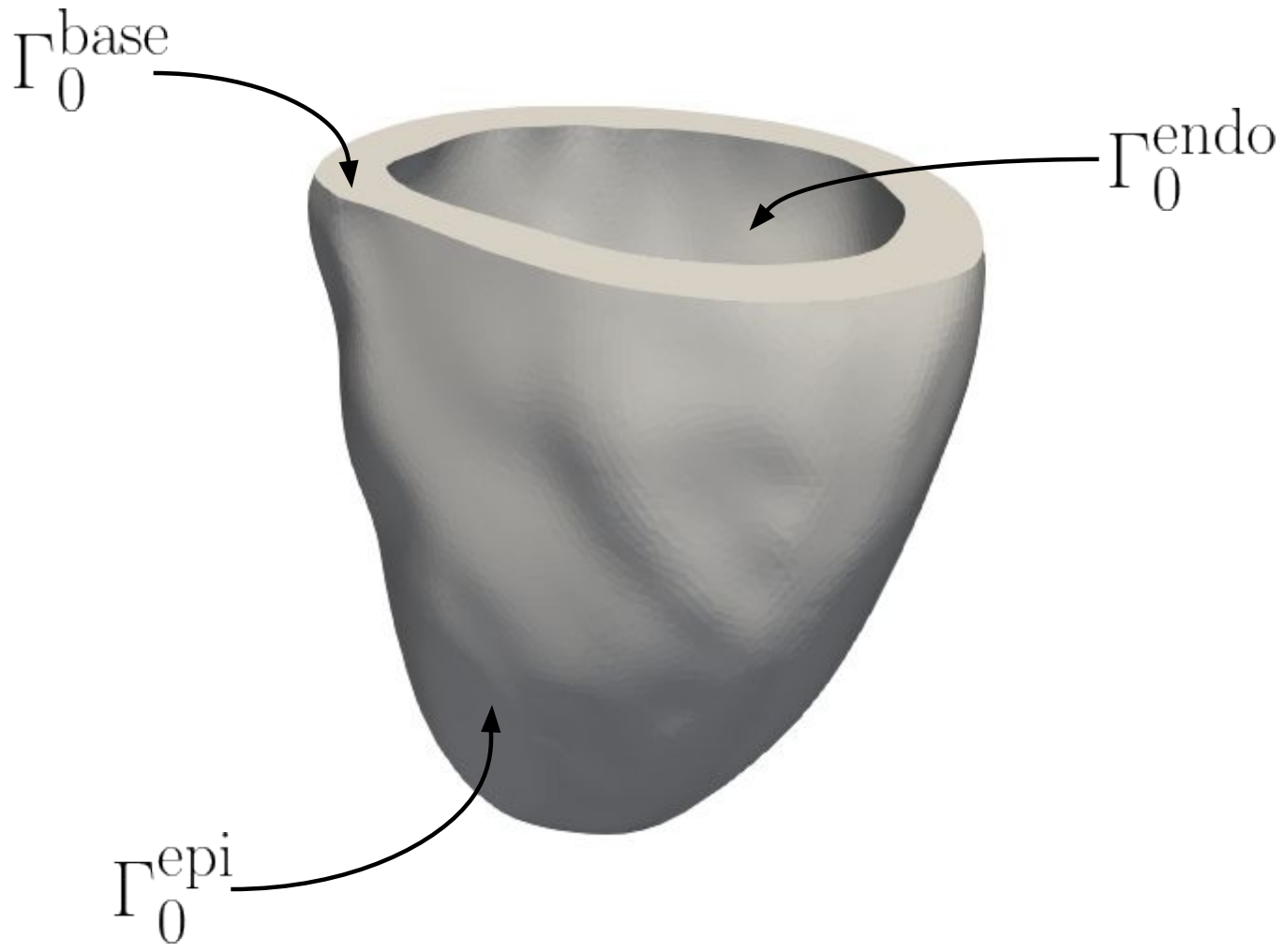

Figure 1: Representation of boundaries $\Gamma_0^{\text{epi}}$, $\Gamma_0^{\text{base}}$ and $\Gamma_0^{\text{endo}}$ of the domain $\Omega_0$, given by the Zygote Solid 3D left ventricle [24].

# 2 Mathematical models

We consider a computational domain $\Omega_0 \subset \mathbb{R}^3$, representing the 3D region occupied by the left ventricle (Fig. 1). We split the boundary of $\Omega_0$ into endocardium ($\Gamma_0^{\text{endo}}$), epicardium ($\Gamma_0^{\text{epi}}$) and ventricular base ($\Gamma_0^{\text{base}}$), namely the artificial boundary located where the left ventricle geometry is cut.

In this paper, we consider a multiphysics and multiscale model of cardiac electromechanics made of five different blocks (henceforth denoted as *core models*) plus a coupling condition. The core models are associated with the different physical phenomena concurring – at different spatial and temporal scales – at the heart function. They correspond to the propagation of the electrical potential ($\mathscr{E}$) [12, 13, 39, 53], ion dynamics ($\mathscr{I}$) [8, 31, 57], active contraction of cardiomyocytes ($\mathscr{A}$) [42, 44, 49–51], tissue mechanics ($\mathscr{M}$) [21, 36] and blood circulation ($\mathscr{C}$) [6, 22]. Finally, the volume conservation condition ($\mathscr{V}$) enables to consistently couple ($\mathscr{M}$) and ($\mathscr{C}$) core models. In Fig. 2 we depict the electric analog circuit corresponding to our 0D circulation model, along with the coupling with a 3D electromechanical description of the left ventricle.

The model features the following unknowns:

$$\begin{aligned} &u\colon \Omega_0 \times (0,T) \to \mathbb{R}, \quad &&\boldsymbol{w}\colon \Omega_0 \times (0,T) \to \mathbb{R}^{n_{\boldsymbol{w}}}, \quad &&\boldsymbol{z}\colon \Omega_0 \times (0,T) \to \mathbb{R}^{n_{\boldsymbol{z}}}, \\ &\mathbf{s}\colon \Omega_0 \times (0,T) \to \mathbb{R}^{n_{\mathbf{s}}}, \quad &&\mathbf{d}\colon \Omega_0 \times (0,T) \to \mathbb{R}^3, \quad &&\boldsymbol{c}_1\colon (0,T) \to \mathbb{R}^{n_{\boldsymbol{c}}}, \\ &p_{\text{LV}}\colon (0,T) \to \mathbb{R} \end{aligned} \tag{1}$$

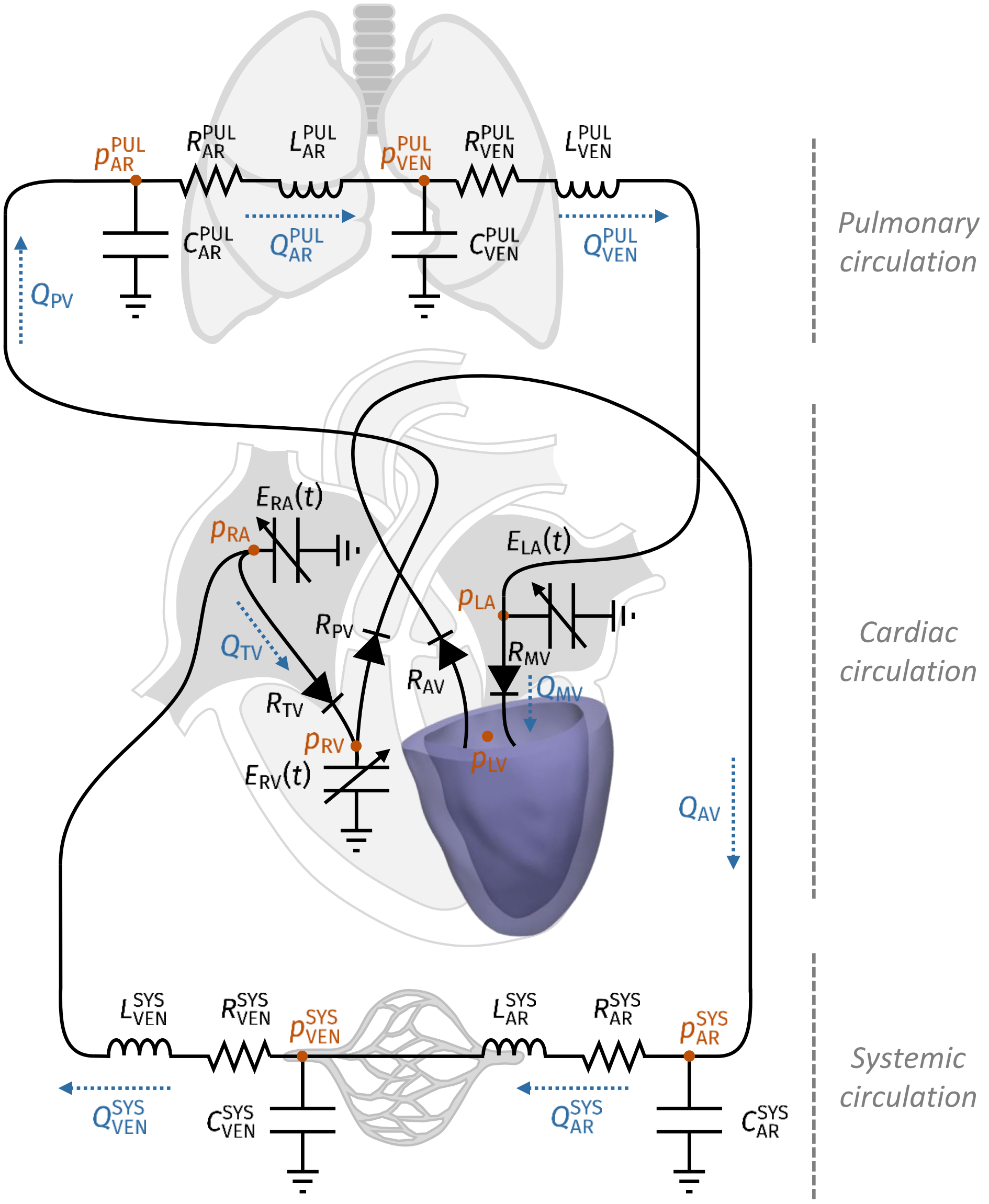

Figure 2: 3D-0D coupling between the left ventricle 3D electromechanical model and the 0D model. The state variables corresponding to pressures and fluxes are depicted in orange and blue, respectively.

where $u$ denotes the transmembrane potential, $\boldsymbol{w}$ and $\boldsymbol{z}$ the ionic variables, $\mathbf{s}$ the states variables of the force generation model, $\mathbf{d}$ the mechanical displacement of the tissue, $\boldsymbol{c}_1$ is the state vector of the circulation models (pressures, volumes and fluxes in the different compartments of the vascular network) and $p_{\mathrm{LV}}$ is the left ventricle pressure. The model reads as follows:

$$
(\mathscr{E})\begin{cases}\chi_{\mathrm{m}}\left[C_{\mathrm{m}}\dfrac{\partial u}{\partial t}+\mathcal{I}_{\mathrm{ion}}(u,\boldsymbol{w},\boldsymbol{z})\right]-\nabla\cdot(J\mathbf{F}^{-1}\boldsymbol{D}_{\mathrm{M}}\mathbf{F}^{-T}\nabla u)\\ \qquad\qquad\qquad\qquad\qquad\qquad\qquad\qquad=\mathcal{I}_{\mathrm{app}}(t),\\ \left(J\mathbf{F}^{-1}\boldsymbol{D}_{\mathrm{M}}\mathbf{F}^{-T}\nabla u\right)\cdot\mathbf{N}=0,\end{cases}\tag{2a}
$$

in $\Omega_0\times(0,T)$, with $u=u_0$ in $\Omega_0$ at time $t=0$.

$$
(\mathscr{I})\begin{cases}\dfrac{\partial\boldsymbol{w}}{\partial t}-\boldsymbol{H}(u,\boldsymbol{w})=\mathbf{0},\\ \dfrac{\partial\boldsymbol{z}}{\partial t}-\boldsymbol{G}(u,\boldsymbol{w},\boldsymbol{z})=\mathbf{0},\end{cases}\tag{2b}
$$

in $\Omega_0\times(0,T)$, with $\boldsymbol{w}=\boldsymbol{w}_0$ and $\boldsymbol{z}=\boldsymbol{z}_0$ in $\Omega_0$ at time $t=0$.

$$
(\mathscr{A})\quad\frac{\partial\mathbf{s}}{\partial t}=\boldsymbol{K}(\mathbf{s},[\mathrm{Ca}^{2+}]_{\mathrm{i}},SL),\tag{2c}
$$

in $\Omega_0\times(0,T)$, with $\mathbf{s}(0)=\mathbf{s}_0$ in $\Omega_0$ at time $t=0$.

$$
(\mathscr{M})\begin{cases}\rho_{\mathrm{s}}\dfrac{\partial^2\mathbf{d}}{\partial t^2}-\nabla\cdot\mathbf{P}(\mathbf{d},T_{\mathrm{a}}(\mathbf{s}))=\mathbf{0},\\ \mathbf{P}(\mathbf{d},T_{\mathrm{a}}(\mathbf{s}))\mathbf{N}+\mathbf{K}^{\mathrm{epi}}\mathbf{d}+\mathbf{C}^{\mathrm{epi}}\dfrac{\partial\mathbf{d}}{\partial t}=\mathbf{0} & \text{on }\Gamma_0^{\mathrm{epi}}\times(0,T),\\ \mathbf{P}(\mathbf{d},T_{\mathrm{a}}(\mathbf{s}))\mathbf{N}=p_{\mathrm{LV}}(t)\,|J\mathbf{F}^{-T}\mathbf{N}|\,\mathbf{v}^{\mathrm{base}}(t) & \text{on }\Gamma_0^{\mathrm{base}}\times(0,T),\\ \mathbf{P}(\mathbf{d},T_{\mathrm{a}}(\mathbf{s}))\mathbf{N}=-p_{\mathrm{LV}}(t)\,J\mathbf{F}^{-T}\mathbf{N} & \text{on }\Gamma_0^{\mathrm{endo}}\times(0,T),\end{cases}\tag{2d}
$$

in $\Omega_0\times(0,T)$, with $\mathbf{d}=\mathbf{d}_0$ and $\dfrac{\partial\mathbf{d}}{\partial t}=\dot{\mathbf{d}}_0$ in $\Omega_0$ at time $t=0$.

$$
(\mathscr{C})\quad\frac{d\boldsymbol{c}_1(t)}{dt}=\widetilde{\boldsymbol{D}}(t,\boldsymbol{c}_1(t),p_{\mathrm{LV}}(t)),\tag{2e}
$$

for $t\in(0,T)$, with $\boldsymbol{c}_1(0)=\boldsymbol{c}_{1,0}$.

$$
(\mathscr{V})\quad V_{\mathrm{LV}}^{\mathrm{0D}}(\boldsymbol{c}_1(t))=V_{\mathrm{LV}}^{\mathrm{3D}}(\mathbf{d}(t)),\tag{2f}
$$

for $t\in(0,T)$.

Hereafter we will denote the pair $(\mathscr{E})$–$(\mathscr{I})$ as *cardiac electrophysiology.* In the following sections, we will detail the different core models, along with the meaning of all the variables, and their coupling.

## 2.1 Electrophysiology $(\mathscr{E})$–$(\mathscr{I})$

We model the propagation of transmembrane potential $u$ with the monodomain equation $(\mathscr{E})$, a diffusion-reaction PDE describing the electric activity of cardiac muscle cells [12, 13, 39, 53]. We couple the monodomain equation with the ten Tusscher-Panfilov ionic model, that we denote by TTP06, due to our focus on the physiological ventricular activity [57]. This model enables to describe the dynamics of ionic fluxes across the cell membrane in an accurate and detailed manner, which is essential for electromechanical coupling [14].

In the electrophysiological model ($\mathscr{E}$)–($\mathscr{I}$), $C_\mathrm{m}$ is the total membrane capacitance and $\chi_\mathrm{m}$ is the area of cell membrane per tissue volume. The vector $\boldsymbol{w} = \{w_1, w_2, ..., w_{n_{\boldsymbol{w}}}\}$ expresses $n_{\boldsymbol{w}}$ recovery (or gating) variables, which play the role of density functions and model the fraction of open ionic channels across the membrane of a single cell. The vector $\boldsymbol{z} = \{z_1, z_2, ..., z_{n_{\boldsymbol{z}}}\}$ defines $n_{\boldsymbol{z}}$ concentration variables of specific ionic species (among which intracellular calcium ions concentration $[\mathrm{Ca}^{2+}]_\mathrm{i}$, which plays a crucial role in the mechanical activation). $\boldsymbol{H}$ and $\boldsymbol{G}$ are suitably defined vector-valued functions, depending on the specific ionic model at hand.

The transmission of the electrical potential $u$ through the so–called gap junctions is described in ($\mathscr{E}$) by means of the diffusion term $\nabla \cdot (J\mathbf{F}^{-1}\boldsymbol{D}_\mathrm{M}\mathbf{F}^{-T}\nabla u)$. In this term, $\boldsymbol{D}_\mathrm{M}$ is the diffusion tensor in the deformed configuration (material coordinates). To account for the dependence of the electrical properties on the tissue on its stretch (mechano-electrical feedback, see e.g. [12, 13]), we introduce the deformation gradient tensor $\mathbf{F} = \mathbf{I} + \nabla \mathbf{d}$, where $\mathbf{d}$ denotes the displacement of the myocardium, and the deformation Jacobian $J = \det(\mathbf{F})$, and we compute the pull-back of $\boldsymbol{D}_\mathrm{M}$ in the reference configuration.

The applied current $\mathcal{I}_\mathrm{app}(t)$ mimics the effect of the Purkinje network [30, 47, 58], which we do not model in this context, by triggering the action potential at specific locations of the myocardium. The ionic current $\mathcal{I}_\mathrm{ion}(u, \boldsymbol{w}, \boldsymbol{z})$ models the multiscale effects from the cellular level up to the tissue one and strictly depends on the selected ionic model. A no-flux (Neumann) boundary condition represents an electrically isolated domain.

To reduce the number of parameters, we rescale the first equation in ($\mathscr{E}$) by $C_\mathrm{m}^{-1}\chi_\mathrm{m}^{-1}$, yielding:

$$\frac{\partial u}{\partial t} + \widetilde{\mathcal{I}}_\mathrm{ion}(u, \boldsymbol{w}, \boldsymbol{z}) = \nabla \cdot (J\mathbf{F}^{-1}\widetilde{\boldsymbol{D}}_\mathrm{M}\mathbf{F}^{-T}\nabla u) + \widetilde{\mathcal{I}}_\mathrm{app}(t), \tag{3}$$

where $\widetilde{\mathcal{I}}_\mathrm{ion} := C_\mathrm{m}^{-1}\mathcal{I}_\mathrm{ion}$, $\widetilde{\boldsymbol{D}}_\mathrm{M} = C_\mathrm{m}^{-1}\chi_\mathrm{m}^{-1}\boldsymbol{D}_\mathrm{M}$ and $\widetilde{\mathcal{I}}_\mathrm{app} := C_\mathrm{m}^{-1}\chi_\mathrm{m}^{-1}\mathcal{I}_\mathrm{app}$. Then, we set $\widetilde{\boldsymbol{D}}_\mathrm{M} = \sigma_\mathrm{t}\mathbf{I} + (\sigma_\mathrm{l} - \sigma_\mathrm{t})\mathbf{f}_0 \otimes \mathbf{f}_0 + (\sigma_\mathrm{n} - \sigma_\mathrm{t})\mathbf{n}_0 \otimes \mathbf{n}_0$ as the diffusion tensor, with $\mathbf{f}_0$ the vector field expressing the fibers direction, $\mathbf{n}_0$ the vector field that indicates the crossfibers direction, and $\sigma_\mathrm{l}, \sigma_\mathrm{t}, \sigma_\mathrm{n} \in \mathbb{R}^+$ longitudinal, transversal and normal conductivities, respectively [12, 52]. $\mathbf{f}_0$, $\mathbf{s}_0$ (i.e. the sheets normal direction) and $\mathbf{n}_0$ account for the anisotropic properties of the cardiac tissue. This local orthonormal coordinate system defined at each point of the computational domain $\Omega_0$ can be generated using rule-based algorithms [4, 17, 38] and has a significant role in the electromechanical framework [3].

## 2.2 Active force generation ($\mathscr{A}$)

Heart contraction is the result of mechano-chemical interactions among the contractile proteins actin and myosin, taking place at the scale of the sarcomeres, the fundamental contractile unit of the cardiac muscle [5, 9, 26, 40, 42]. To model such complex mechanisms, we consider the model that we proposed in [41], henceforth denoted as RDQ18 model. This model is based on a biophysically detailed description of the sarcomeric proteins with an explicit representation of the cooperative nearest-neighboring interactions, responsible for the high sensitivity of the cardiac contractile apparatus to changes in calcium concentration, which is one of the ionic species modelled by the TTP06 ionic model. Thanks to the spatially-explicit representation of the sarcomere filaments, the RDQ18 model also incorporates the feedback of the sarcomere shortening, resulting from the muscle contraction, on the force generation mechanism itself. This is of outmost importance since the regulation due to the sarcomere length ($SL$) constitutes the microscopical basis of the well-known Frank-Starling mechanism at the macroscale level; in practice, higher end-diastolic volumes translate into higher stroke volumes [26].

The RDQ18 model is a system of ODEs in the form of ($\mathscr{A}$), where the vector $\mathbf{s}$ collects the variables of the RDQ18 model and where $\boldsymbol{K}$ is a suitable function defined in [41] (we remark that $\boldsymbol{K}$ does not involve derivatives of $\mathbf{s}$ with respect to the spatial variable). Within a multiscale framework, the RDQ18 model is ideally set at every point of the computational domain $\Omega_0$. The input variable $[\mathrm{Ca}^{2+}]_\mathrm{i}$ is provided by the TTP06 ionic model in each point of the domain, while $SL$ is provided by the solution of the mechanical model, as we will explain later in Sec. 2.3.

The output of the RDQ18 model is the permissivity $P \in [0, 1]$, obtained as a function of the states $\mathbf{s}$ (i.e. $P = G(\mathbf{s})$, where $G$ is a linear function defined in [41]). The permissivity represents the fraction of contractile units being in the force-generating state. Hence, the effective active tension is given by $T_\mathrm{a} = T_\mathrm{a}^{\max} P$, where $T_\mathrm{a}^{\max}$ denotes the tension generated when all the contractile units are generating force (i.e. for $P = 1$).

The RDQ18 model accurately describes the microscopic force generation mechanisms. This accuracy results in a higher computational cost compared to phenomenological models typically used for multiscale simulations (see e.g. [29, 33]). To overcome this issue, in the multiscale model of electromechanics we take advantage of the model based on Artificial Neural Networks (ANNs) presented in [44]. This model is a fast surrogate of the RDQ18 model (high-fidelity model), learned from a collection of pre-computed simulations obtained with the RDQ18 model itself, thanks to the Machine Learning algorithm that we proposed in [43]. Such reduced model is written in the same form of ($\mathscr{A}$). However, the state vector $\mathbf{s}$ of ANN-based model only contains two variables, instead of the more than 2000 variables of the high fidelity model. This significantly reduces the computational costs associated with its numerical approximation (both in terms of CPU time and memory storage), at the price of only a small approximation, as the overall relative error between the results of the high-fidelity and of the reduced models is of the order of $10^{-3}$ [44]. In this way we obtain an excellent trade-off between computational cost and biophysical accuracy of the results.

### 2.3 Active and passive mechanics ($\mathscr{M}$)

We describe the dynamics of the tissue displacement $\mathbf{d}$ by the momentum conservation equation ($\mathscr{M}$) (see e.g. [36]). The Piola-Kirchhoff stress tensor $\mathbf{P} = \mathbf{P}(\mathbf{d}, T_\mathrm{a})$ incorporates both passive and active mechanics of the tissue. Under the hyperelasticity assumption, once the strain energy density function $\mathcal{W} : \mathrm{Lin}^+ \to \mathbb{R}$ is introduced, the passive part of the Piola-Kirchhoff stress tensor is obtained as $\frac{\partial}{\partial \mathbf{F}} \mathcal{W}(\mathbf{F})$. In conclusion, the full Piola-Kirchhoff tensor reads:

$$\mathbf{P}(\mathbf{d}, T_\mathrm{a}) = \frac{\partial \mathcal{W}(\mathbf{F})}{\partial \mathbf{F}} + T_\mathrm{a} \frac{\mathbf{F}\mathbf{f}_0 \otimes \mathbf{f}_0}{\sqrt{\mathcal{I}_{4f}}}, \tag{4}$$

where the first term stands as the passive part, while the latter as the active one of the tensor $\mathbf{P}$, and where $T_\mathrm{a}$ denotes the active tension, provided the force generation model of Sec. 2.2. $\mathcal{I}_{4f} = \mathbf{F}\mathbf{f}_0 \cdot \mathbf{F}\mathbf{f}_0$ is a measure of the tissue stretch along the fibers direction.

Several models have been proposed in literature to describe the anisotropic nature of the cardiac muscle tissue. In this paper, we consider the Guccione strain energy density function [20, 21], that reads $\mathcal{W}(\mathbf{F}) = \frac{C}{2}\left(e^Q - 1\right)$, with

$$Q = b_\mathrm{ff} E_\mathrm{ff}^2 + b_\mathrm{ss} E_\mathrm{ss}^2 + b_\mathrm{nn} E_\mathrm{nn}^2 + b_\mathrm{fs}\left(E_\mathrm{fs}^2 + E_\mathrm{sf}^2\right) + b_\mathrm{fn}\left(E_\mathrm{fn}^2 + E_\mathrm{nf}^2\right) + b_\mathrm{sn}\left(E_\mathrm{sn}^2 + E_\mathrm{ns}^2\right), \tag{5}$$

where $E_\mathrm{ab} = \boldsymbol{E}\boldsymbol{a}_0 \cdot \boldsymbol{b}_0 \quad \text{for} \quad a, b \in \{f, s, n\}$ are the entries of $\mathbf{E} = \frac{1}{2}(\mathbf{C} - \mathbf{I})$, i.e the Green-Lagrange strain energy tensor, being $\mathbf{C} = \mathbf{F}^T\mathbf{F}$ the right Cauchy-Green deformation tensor. We consider a further term, defined as $\mathcal{W}_\mathrm{vol}(J) = \frac{B}{2}(J - 1)\log(J)$, convex in $J$ and such that $J = 1$ is the global minimum, which penalizes large variations of volume, thus realizing a (weakly) incompressible constraint [11, 16, 61]; $B \in \mathbb{R}^+$ represents the bulk modulus.

To model the interaction of the left ventricle with the pericardium [19, 37, 56], we impose at the epicardial boundary $\Gamma_0^{\text{epi}}$ the generalized Robin boundary condition $\mathbf{PN} + \mathbf{K}^{\text{epi}}\mathbf{d} + \mathbf{C}^{\text{epi}}\frac{\partial \mathbf{d}}{\partial t} = \mathbf{0}$, by defining the following tensors

$$\mathbf{K}^{\text{epi}} = K_{\perp}^{\text{epi}}(\mathbf{N} \otimes \mathbf{N}) + K_{\parallel}^{\text{epi}}(\mathbf{I} - \mathbf{N} \otimes \mathbf{N}),$$
$$\mathbf{C}^{\text{epi}} = C_{\perp}^{\text{epi}}(\mathbf{N} \otimes \mathbf{N}) + C_{\parallel}^{\text{epi}}(\mathbf{I} - \mathbf{N} \otimes \mathbf{N}),$$

where the constants $K_{\perp}^{\text{epi}}$, $K_{\parallel}^{\text{epi}}$, $C_{\perp}^{\text{epi}}$, $C_{\parallel}^{\text{epi}} \in \mathbb{R}^+$ are local values of stiffness and viscosity of the epicardial tissue in the normal or tangential directions, respectively. At the base $\Gamma_0^{\text{base}}$, we set the energy-consistent boundary condition $\mathbf{PN} = p_{\text{LV}}(t)\,|J\mathbf{F}^{-T}\mathbf{N}|\,\mathbf{v}^{\text{base}}(t)$, originally proposed in [44], that provides an explicit expression for the stresses located at the artificial boundary $\Gamma_0^{\text{base}}$, where we have defined the vector

$$\mathbf{v}^{\text{base}}(t) = \frac{\int_{\Gamma_0^{\text{endo}}} J\mathbf{F}^{-T}\mathbf{N} d\Gamma_0}{\int_{\Gamma_0^{\text{base}}} |J\mathbf{F}^{-T}\mathbf{N}| d\Gamma_0}.$$

As we later show (Sec. 4.2), this formulation allows to straightforwardly couple the 3D mechanical model with a 0D model of the whole circulation in an energetically consistent manner. Finally, at the endocardium $\Gamma_0^{\text{endo}}$, the boundary condition $\mathbf{PN} = -p_{\text{LV}}(t)\, J\mathbf{F}^{-T}\mathbf{N}$ accounts for the pressure $p_{\text{LV}}(t)$ exerted by the blood contained in the ventricular chamber, modeled through 0D closed-loop circulation model.

As anticipated, the mechanical model ($\mathscr{M}$) has a feedback on the force generation model of ($\mathscr{A}$), as $\mathbf{d}$ determines the local sarcomere length $SL$. More precisely, since sarcomeres are aligned with the muscle fibers $\mathbf{f}_0$, the local sarcomere length $SL$ is given as $SL = SL_0\sqrt{\mathcal{I}_{4f}}$, where $SL_0$ denotes the sarcomere length at rest. To recover the $SL$ field, we consider the following differential problem:

$$\begin{cases} -\delta_{SL}^2 \Delta SL + SL = SL_0\sqrt{\mathcal{I}_{4f}} & \text{in } \Omega_0 \times (0,T), \\ \delta_{SL}^2 \nabla SL \cdot \mathbf{N} = 0 & \text{on } \partial\Omega_0 \times (0,T), \end{cases} \tag{6}$$

where $\delta_{SL}$ is the regularization parameter, whose aim is that of making smoother the field $SL$ across a computational domain, preventing sharp spatial variations across scales smaller than $\delta_{SL}$. This will be particularly useful in view of its FEM approximation in Part II [45].

## 2.4 Blood circulation ($\mathscr{C}$)

To model the hemodynamics of the whole circulatory network, we propose a lumped parameters closed-loop model, inspired by previous models available in literature [6, 22]. Systemic and pulmonary circulations are modeled with resistance-inductance-capacitance (RLC) circuits, one for the arterial part and the other one for the venous part. The four chambers are modeled by time-varying elastance elements, whereas the four valves are represented as non-ideal diodes. Our 0D closed-loop

circulation model reads:

$$
\begin{cases}
\dfrac{dV_{\mathrm{LA}}(t)}{dt} = Q_{\mathrm{VEN}}^{\mathrm{PUL}}(t) - Q_{\mathrm{MV}}(t), \qquad \dfrac{dV_{\mathrm{LV}}(t)}{dt} = Q_{\mathrm{MV}}(t) - Q_{\mathrm{AV}}(t), & \text{(7a)}\\
\dfrac{dV_{\mathrm{RA}}(t)}{dt} = Q_{\mathrm{VEN}}^{\mathrm{SYS}}(t) - Q_{\mathrm{TV}}(t), \qquad \dfrac{dV_{\mathrm{RV}}(t)}{dt} = Q_{\mathrm{TV}}(t) - Q_{\mathrm{PV}}(t), & \text{(7b)}\\
C_{\mathrm{AR}}^{\mathrm{SYS}} \dfrac{dp_{\mathrm{AR}}^{\mathrm{SYS}}(t)}{dt} = Q_{\mathrm{AV}}(t) - Q_{\mathrm{AR}}^{\mathrm{SYS}}(t), \quad C_{\mathrm{VEN}}^{\mathrm{SYS}} \dfrac{dp_{\mathrm{VEN}}^{\mathrm{SYS}}(t)}{dt} = Q_{\mathrm{AR}}^{\mathrm{SYS}}(t) - Q_{\mathrm{VEN}}^{\mathrm{SYS}}(t), & \text{(7c)}\\
C_{\mathrm{AR}}^{\mathrm{PUL}} \dfrac{dp_{\mathrm{AR}}^{\mathrm{PUL}}(t)}{dt} = Q_{\mathrm{PV}}(t) - Q_{\mathrm{AR}}^{\mathrm{PUL}}(t), \quad C_{\mathrm{VEN}}^{\mathrm{PUL}} \dfrac{dp_{\mathrm{VEN}}^{\mathrm{PUL}}(t)}{dt} = Q_{\mathrm{AR}}^{\mathrm{PUL}}(t) - Q_{\mathrm{VEN}}^{\mathrm{PUL}}(t), & \text{(7d)}\\
\dfrac{L_{\mathrm{AR}}^{\mathrm{SYS}}}{R_{\mathrm{AR}}^{\mathrm{SYS}}} \dfrac{dQ_{\mathrm{AR}}^{\mathrm{SYS}}(t)}{dt} = -Q_{\mathrm{AR}}^{\mathrm{SYS}}(t) - \dfrac{p_{\mathrm{VEN}}^{\mathrm{SYS}}(t) - p_{\mathrm{AR}}^{\mathrm{SYS}}(t)}{R_{\mathrm{AR}}^{\mathrm{SYS}}}, & \text{(7e)}\\
\dfrac{L_{\mathrm{VEN}}^{\mathrm{SYS}}}{R_{\mathrm{VEN}}^{\mathrm{SYS}}} \dfrac{dQ_{\mathrm{VEN}}^{\mathrm{SYS}}(t)}{dt} = -Q_{\mathrm{VEN}}^{\mathrm{SYS}}(t) - \dfrac{p_{\mathrm{RA}}(t) - p_{\mathrm{VEN}}^{\mathrm{SYS}}(t)}{R_{\mathrm{VEN}}^{\mathrm{SYS}}}, & \text{(7f)}\\
\dfrac{L_{\mathrm{AR}}^{\mathrm{PUL}}}{R_{\mathrm{AR}}^{\mathrm{PUL}}} \dfrac{dQ_{\mathrm{AR}}^{\mathrm{PUL}}(t)}{dt} = -Q_{\mathrm{AR}}^{\mathrm{PUL}}(t) - \dfrac{p_{\mathrm{VEN}}^{\mathrm{PUL}}(t) - p_{\mathrm{AR}}^{\mathrm{PUL}}(t)}{R_{\mathrm{AR}}^{\mathrm{PUL}}}, & \text{(7g)}\\
\dfrac{L_{\mathrm{VEN}}^{\mathrm{PUL}}}{R_{\mathrm{VEN}}^{\mathrm{PUL}}} \dfrac{dQ_{\mathrm{VEN}}^{\mathrm{PUL}}(t)}{dt} = -Q_{\mathrm{VEN}}^{\mathrm{PUL}}(t) - \dfrac{p_{\mathrm{LA}}(t) - p_{\mathrm{VEN}}^{\mathrm{PUL}}(t)}{R_{\mathrm{VEN}}^{\mathrm{PUL}}}, & \text{(7h)}
\end{cases}
$$

with $t \in (0, T)$, where:

$$p_{\mathrm{LV}}(t) = p_{\mathrm{EX}}(t) + E_{\mathrm{LV}}(t) \left(V_{\mathrm{LV}}(t) - V_{0,\mathrm{LV}}\right), \tag{8a}$$

$$p_{\mathrm{LA}}(t) = p_{\mathrm{EX}}(t) + E_{\mathrm{LA}}(t) \left(V_{\mathrm{LA}}(t) - V_{0,\mathrm{LA}}\right), \tag{8b}$$

$$p_{\mathrm{RV}}(t) = p_{\mathrm{EX}}(t) + E_{\mathrm{RV}}(t) \left(V_{\mathrm{RV}}(t) - V_{0,\mathrm{RV}}\right), \tag{8c}$$

$$p_{\mathrm{RA}}(t) = p_{\mathrm{EX}}(t) + E_{\mathrm{RA}}(t) \left(V_{\mathrm{RA}}(t) - V_{0,\mathrm{RA}}\right), \tag{8d}$$

$$Q_{\mathrm{MV}}(t) = \frac{p_{\mathrm{LA}}(t) - p_{\mathrm{LV}}(t)}{R_{\mathrm{MV}}(p_{\mathrm{LA}}(t), p_{\mathrm{LV}}(t))}, \qquad Q_{\mathrm{AV}}(t) = \frac{p_{\mathrm{LV}}(t) - p_{\mathrm{AR}}^{\mathrm{SYS}}(t)}{R_{\mathrm{AV}}(p_{\mathrm{LV}}(t), p_{\mathrm{AR}}^{\mathrm{SYS}}(t))}, \tag{8e}$$

$$Q_{\mathrm{TV}}(t) = \frac{p_{\mathrm{RA}}(t) - p_{\mathrm{RV}}(t)}{R_{\mathrm{TV}}(p_{\mathrm{RA}}(t), p_{\mathrm{RV}}(t))}, \qquad Q_{\mathrm{PV}}(t) = \frac{p_{\mathrm{RV}}(t) - p_{\mathrm{AR}}^{\mathrm{PUL}}(t)}{R_{\mathrm{PV}}(p_{\mathrm{RV}}(t), p_{\mathrm{AR}}^{\mathrm{PUL}}(t))}, \tag{8f}$$

with $t \in (0, T)$. In this model, $p_{\mathrm{LA}}(t)$, $p_{\mathrm{RA}}(t)$, $p_{\mathrm{LV}}(t)$, $p_{\mathrm{RV}}(t)$, $V_{\mathrm{LA}}(t)$, $V_{\mathrm{RA}}(t)$, $V_{\mathrm{LV}}(t)$ and $V_{\mathrm{RV}}(t)$ refer to pressures and volumes in left atrium, right atrium, left ventricle and right ventricle, respectively. The variables $Q_{\mathrm{MV}}(t)$, $Q_{\mathrm{AV}}(t)$, $Q_{\mathrm{TV}}(t)$ and $Q_{\mathrm{PV}}(t)$ indicate the flow rates through mitral, aortic, tricuspid and pulmonary valves, respectively. Moreover, $p_{\mathrm{AR}}^{\mathrm{SYS}}(t)$, $Q_{\mathrm{AR}}^{\mathrm{SYS}}(t)$, $p_{\mathrm{VEN}}^{\mathrm{SYS}}(t)$ and $Q_{\mathrm{VEN}}^{\mathrm{SYS}}(t)$ express pressures and flow rates of the systemic circulation (arterial and venous). Similarly, $p_{\mathrm{AR}}^{\mathrm{PUL}}(t)$, $Q_{\mathrm{AR}}^{\mathrm{PUL}}(t)$, $p_{\mathrm{VEN}}^{\mathrm{PUL}}(t)$ and $Q_{\mathrm{VEN}}^{\mathrm{PUL}}(t)$ define pressures and flow rates of the pulmonary circulation (arterial and venous). $p_{\mathrm{EX}}(t)$ represents the pressure exerted outside the heart by the surrounding organs and respiration. Time varying $E_{\mathrm{LA}}(t)$, $E_{\mathrm{LV}}(t)$, $E_{\mathrm{RA}}(t)$, $E_{\mathrm{RV}}(t)$ are the analytically prescribed elastances of the four cardiac chambers calibrated on a physiological basis, with values ranging from $E_{\mathrm{LA}}^{\mathrm{pass}}$, $E_{\mathrm{LV}}^{\mathrm{pass}}$, $E_{\mathrm{RA}}^{\mathrm{pass}}$, $E_{\mathrm{RV}}^{\mathrm{pass}}$ – when the chambers are at rest – to $(E_{\mathrm{LA}}^{\mathrm{pass}} + E_{\mathrm{LA}}^{\mathrm{act,max}})$, $(E_{\mathrm{LV}}^{\mathrm{pass}} + E_{\mathrm{LV}}^{\mathrm{act,max}})$, $(E_{\mathrm{RA}}^{\mathrm{pass}} + E_{\mathrm{RA}}^{\mathrm{act,max}})$, $(E_{\mathrm{RV}}^{\mathrm{pass}} + E_{\mathrm{RV}}^{\mathrm{act,max}})$ – when the chambers are fully contracted. Finally, $R_{\mathrm{MV}}(p_1, p_2)$, $R_{\mathrm{AV}}(p_1, p_2)$, $R_{\mathrm{TV}}(p_1, p_2)$ and $R_{\mathrm{PV}}(p_1, p_2)$ define the behavior of valves as diodes, according to the following relationship:

$$R_{\mathrm{i}}(p_1, p_2) = \begin{cases} R_{\min}, & p_1 < p_2 \\ R_{\max}, & p_1 \geq p_2 \end{cases} \quad \text{for} \quad i \in \{\mathrm{MV}, \mathrm{AV}, \mathrm{TV}, \mathrm{PV}\},$$

where $p_1$ and $p_2$ denote the pressures ahead and behind the valve leaflets with respect to the flow direction, whereas $R_{\min}$ and $R_{\max}$ are the minimum and maximum resistance of the valves. For an idealized valve, one would have $R_{\min} = 0$ and $R_{\max} = +\infty$ instead. By setting $R_{\min} > 0$, one has dissipation of mechanical energy taking place when the blood flows through the opened valve (see Sec. 4); we set $R_{\max} < +\infty$ sufficiently large so that blood leakage when the valve is closed is negligible.

Hereafter, for the sake of brevity, Eqs. (7)–(8) will be expressed in the following form:

$$\begin{cases} \dfrac{d\boldsymbol{c}_1(t)}{dt} = \boldsymbol{D}(t, \boldsymbol{c}_1(t), \boldsymbol{c}_2(t)) & t \in (0, T], \\ \boldsymbol{c}_2(t) = \boldsymbol{W}(t, \boldsymbol{c}_1(t)) & t \in [0, T], \\ \boldsymbol{c}_1(0) = \boldsymbol{c}_{1,0}, \end{cases} \tag{9}$$

where:

$$\begin{aligned} \boldsymbol{c}_1(t) = {} & (V_{\mathrm{LA}}(t), V_{\mathrm{LV}}(t), V_{\mathrm{RA}}(t), V_{\mathrm{RV}}(t), p_{\mathrm{AR}}^{\mathrm{SYS}}(t), p_{\mathrm{VEN}}^{\mathrm{SYS}}(t), p_{\mathrm{AR}}^{\mathrm{PUL}}(t), p_{\mathrm{VEN}}^{\mathrm{PUL}}(t), \\ & Q_{\mathrm{AR}}^{\mathrm{SYS}}(t), Q_{\mathrm{VEN}}^{\mathrm{SYS}}(t), Q_{\mathrm{AR}}^{\mathrm{PUL}}(t), Q_{\mathrm{VEN}}^{\mathrm{PUL}}(t))^T, \\ \boldsymbol{c}_2(t) = {} & (p_{\mathrm{LV}}(t), p_{\mathrm{LA}}(t), p_{\mathrm{RV}}(t), p_{\mathrm{RA}}(t), Q_{\mathrm{MV}}(t), Q_{\mathrm{AV}}(t), Q_{\mathrm{TV}}(t), Q_{\mathrm{PV}}(t))^T; \end{aligned}$$

$\boldsymbol{D}(t, \boldsymbol{c}_1(t), \boldsymbol{c}_2(t))$ collects the whole r.h.s. of Eq. (7), while $\boldsymbol{c}_2(t) = \boldsymbol{W}(t, \boldsymbol{c}_1(t))$ stands as a compact notation for Eq. (8), rewritten in explicit form with respect to the variable $\boldsymbol{c}_2$.

## 2.5 3D-0D coupling $(\mathscr{V})$

In Eq. (9) each cardiac chamber is modeled as a time-varying elastance element, that is a 0D simplified model. In this paper, we employ this 0D circulation model in conjunction with the 3D left ventricular model given by $(\mathscr{E})$–$(\mathscr{I})$–$(\mathscr{A})$–$(\mathscr{M})$. With this goal, we remove from the circulation model the time-varying elastance element associated with the left ventricle, and we replace it with the 3D electromechanical model. Hence, the pressure-volume relationship between $p_{\mathrm{LV}}$ and $V_{\mathrm{LV}}$ is no longer prescribed by Eq. (8a), but by the resolution of the 3D electromechanical model. The resulting 3D–0D coupled model (depicted in Fig. 2) must satisfy at each time $t \in (0, T)$ the volume-consistency coupling condition $V_{\mathrm{LV}}^{\mathrm{0D}}(\boldsymbol{c}_1(t)) = V_{\mathrm{LV}}^{\mathrm{3D}}(\mathbf{d}(t))$, that we denote by $(\mathscr{V})$, where $V_{\mathrm{LV}}^{\mathrm{0D}}(\boldsymbol{c}_1(t)) = V_{\mathrm{LV}}(t)$ represents the left ventricle volume in the 0D circulation model. $V_{\mathrm{LV}}^{\mathrm{3D}}(\mathbf{d}(t))$ represents the left ventricle volume in the 3D model and it is computed as:

$$V_{\mathrm{LV}}^{\mathrm{3D}}(\mathbf{d}(t)) = \int_{\Gamma_0^{\mathrm{endo}}} J(t) \left( (\mathbf{h} \otimes \mathbf{h}) \left( \mathbf{x} + \mathbf{d}(t) - \mathbf{b} \right) \right) \cdot \mathbf{F}^{-T}(t) \mathbf{N} \, d\Gamma_0,$$

where $\mathbf{h}$ is a vector orthogonal to the left ventricle centerline (i.e. lying on the left ventricle base) [48]. Subtracting to the space coordinate $\mathbf{x} + \mathbf{d}(t)$ that of a point $\mathbf{b}$, lying inside the left ventricle, improves the accuracy of the formula when the ventricular base changes its orientation.

Having introduced an additional scalar equation, i.e. $(\mathscr{V})$, we expect an additional unknown: it is in fact $p_{\mathrm{LV}}$, which is not determined by Eq. (8a) anymore. Rather, it acts as a Lagrange multiplier enforcing the constraint $(\mathscr{V})$.

Hence, we define the "reduced" vector $\widetilde{\boldsymbol{c}}_2$ such that $\boldsymbol{c}_2^T = (p_{\mathrm{LV}}, \widetilde{\boldsymbol{c}}_2^T)$, so that we can rewrite Eqs. (8a)–(8f) as $\widetilde{\boldsymbol{c}}_2(t) = \widetilde{\boldsymbol{W}}(t, \boldsymbol{c}_1(t), p_{\mathrm{LV}}(t))$. This allows to write the "reduced" version of Eq. (9) as $(\mathscr{C})$ where we have defined

$$\widetilde{\boldsymbol{D}}(t, \boldsymbol{c}_1, p_{\mathrm{LV}}) := \boldsymbol{D}\left(t, \boldsymbol{c}_1, \begin{pmatrix} p_{\mathrm{LV}} \\ \widetilde{\boldsymbol{W}}(t, \boldsymbol{c}_1, p_{\mathrm{LV}}) \end{pmatrix}\right).$$

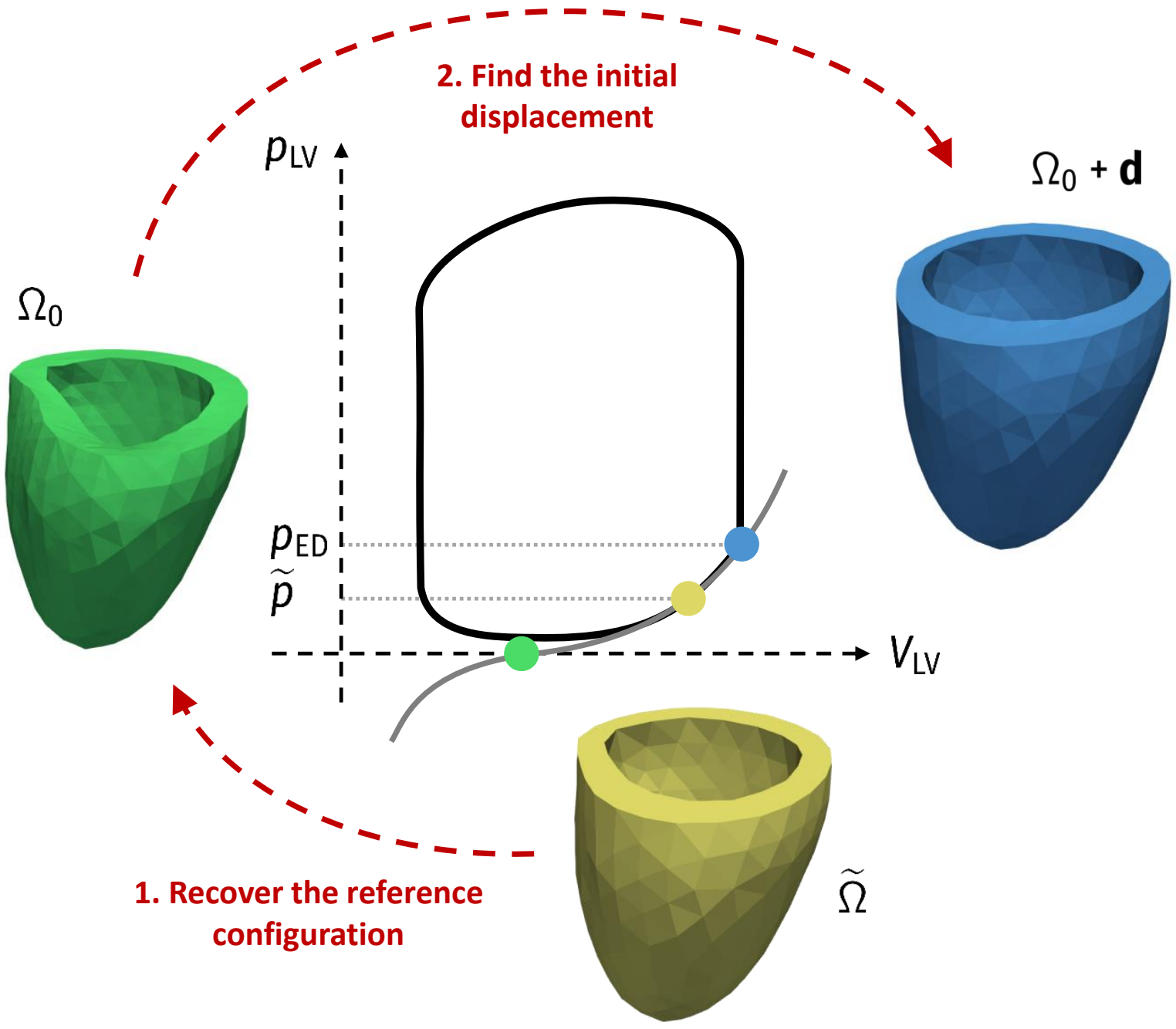

Figure 3: Sketch of the strategy used to initialize the simulation. The grey line represents the so–called Klotz curve [27], that is the pressure-volume relationship of the relaxed ventricle. The black line represents the pressure-volume loop of the left ventricle.

In conclusion, we obtain the coupled 3D-0D model reported in Eq. (2). We remark that the number of equations balances with the number of unknowns defined in (1): we have $1+n_{\boldsymbol{w}}+n_{\boldsymbol{z}}+n_{\mathbf{s}}+3$ unknowns (respectively, equations) defined in $\Omega_0 \times (0,T)$ and $n_{\mathbf{c}}+1$ unknowns (respectively, equations) defined in $(0,T)$.

# 3 Reference configuration and initial displacement

In the mechanical model ($\mathscr{M}$), the stress-strain relationship (4) is referred to the natural stress-free configuration $\Omega_0$. However, the geometry is never unloaded during the cardiac cycle. Take for example the case of medical images to generate the left ventricle: this geometry, that we denote by $\widetilde{\Omega}$, does not correspond to the stress-free configuration, since an internal pressure $p_{\mathrm{LV}} \neq 0$ occurs at every stage of the heartbeat. Therefore, in the preprocessing stage, we need to recover the reference configuration $\Omega_0$ from $\widetilde{\Omega}$.

Our strategy to initialize the numerical simulation in a physically sound manner is sketched in Fig. 3. As a first step, starting from the geometry acquired from medical imaging $\widetilde{\Omega}$, we recover the stress-free reference configuration $\Omega_0$, by virtually deflating the left ventricle, previously subject to an internal pressure $\widetilde{p}$. Then, as a second step, we inflate the left ventricle again, by applying the end diastolic pressure $p_{\mathrm{ED}}$ at the endocardium. In the next sections we give the mathematical details of these two steps.

### 3.1 Recovering the reference configuration

We assume that the configuration $\widetilde{\Omega}$ occurs during diastole, when the left ventricle is loaded by a small value of the pressure $p_{\mathrm{LV}} = \widetilde{p}$ and only a residual active tension $T_{\mathrm{a}} = \widetilde{T}_{\mathrm{a}} > 0$ acts. By adopting a quasi-static assumption (motivated by the slow movement of the myocardium during the final part of diastole), the tissue displacement is given by the solution of the following differential problem:

$$\begin{cases} \nabla \cdot \mathbf{P}(\mathbf{d}, T_{\mathrm{a}}) = \mathbf{0} & \text{in } \Omega_0, \\ \mathbf{P}(\mathbf{d}, T_{\mathrm{a}})\mathbf{N} + \mathbf{K}^{\mathrm{epi}}\mathbf{d} = \mathbf{0} & \text{on } \Gamma_0^{\mathrm{epi}}, \\ \mathbf{P}(\mathbf{d}, T_{\mathrm{a}})\mathbf{N} = p_{\mathrm{LV}}\, |J\mathbf{F}^{-T}\mathbf{N}|\, \mathbf{v}^{\mathrm{base}}(t) & \text{on } \Gamma_0^{\mathrm{base}}, \\ \mathbf{P}(\mathbf{d}, T_{\mathrm{a}})\mathbf{N} = -p_{\mathrm{LV}}\, J\mathbf{F}^{-T}\mathbf{N} & \text{on } \Gamma_0^{\mathrm{endo}}, \end{cases} \tag{10}$$

derived from $(\mathscr{M})$ by setting to zero the time-dependent terms. Thus, to recover the coordinate $\mathbf{x}_0$ of the configuration $\Omega_0$ we need to solve the following inverse problem: find the domain $\Omega_0$ such that, if we displace $\mathbf{x}_0$ by the solution $\mathbf{d}$ of Eq. (10) obtained for $p_{\mathrm{LV}} = \widetilde{p}$ and $T_{\mathrm{a}} = \widetilde{T}_{\mathrm{a}}$, we get the coordinate $\widetilde{\mathbf{x}}$ of the domain $\widetilde{\Omega}$ (i.e. $\widetilde{\mathbf{x}} = \mathbf{x}_0 + \mathbf{d}$). In Part II of this paper [45] we present an algorithm for its numerical solution.

### 3.2 Finding the initial displacement

After the recovery of the reference configuration $\Omega_0$, we set $p_{\mathrm{LV}}$ equal to the end diastolic pressure and we solve again Eq. (10). In this manner, we obtain the end-diastolic configuration of the left ventricle. Hence, we employ the solution $\mathbf{d}$ as initial condition $\mathbf{d}_0$ for $\mathbf{d}$ in $(\mathscr{M})$.

## 4 On the balance of mechanical energy of the electromechanical model

We derive a balance of the mechanical energy of the closed-loop circulation model (9) and we highlight energy injection, dissipation and transfer in the different compartments and in the different stages of the heartbeat. First, in Sec. 4.1, we consider the 0D circulation model introduced in Sec. 2.4. Then, in Sec. 4.2 we consider the coupled 3D-0D model, showing that our formulation is compliant with the above mentioned balance of the mechanical energy.

### 4.1 Energy balance for the 0D model

To define the terms associated with the work performed by the cardiac chambers, we write $E_{\mathrm{i}}(t) = E_{\mathrm{i}}^{\mathrm{pass}} + E_{\mathrm{i}}^{\mathrm{act}}(t)$ (for $i \in \{\mathrm{LA}, \mathrm{LV}, \mathrm{RA}, \mathrm{RV}\}$), where $E_{\mathrm{i}}^{\mathrm{pass}}$ is the passive elastance of the tissue (i.e. the elastance when the tissue is not activated) and $E_{\mathrm{i}}^{\mathrm{act}}$ is instead the active component of the elastance.

**Definition 1.** *We define the total mechanical energy of the whole 0D circulation model as*

$$\begin{aligned} \mathcal{M}(t) = {} & \mathcal{E}_{\mathrm{LA}}(t) + \mathcal{E}_{\mathrm{LV}}(t) + \mathcal{E}_{\mathrm{RA}}(t) + \mathcal{E}_{\mathrm{RV}}(t) + \mathcal{E}_{\mathrm{AR}}^{\mathrm{SYS}}(t) + \mathcal{E}_{\mathrm{VEN}}^{\mathrm{SYS}}(t) + \mathcal{E}_{\mathrm{AR}}^{\mathrm{PUL}}(t) + \mathcal{E}_{\mathrm{VEN}}^{\mathrm{PUL}}(t) \\ & + \mathcal{K}_{\mathrm{AR}}^{\mathrm{SYS}}(t) + \mathcal{K}_{\mathrm{VEN}}^{\mathrm{SYS}}(t) + \mathcal{K}_{\mathrm{AR}}^{\mathrm{PUL}}(t) + \mathcal{K}_{\mathrm{VEN}}^{\mathrm{PUL}}(t), \end{aligned}$$

*where, for* $i \in \{\mathrm{LA}, \mathrm{LV}, \mathrm{RA}, \mathrm{RV}\}$, $j \in \{\mathrm{AR}, \mathrm{VEN}\}$ *and* $k \in \{\mathrm{SYS}, \mathrm{PUL}\}$*:*

- $\mathcal{E}_i(t) = \dfrac{1}{2} E_i^{\mathrm{pass}} \left(V_i(t) - V_{0,i}\right)^2$ *is the elastic energy stored by a cardiac chamber;*

- $\mathcal{E}_j^k(t) = \dfrac{1}{2} C_j^k \left(p_j^k(t)\right)^2$ *is the elastic energy stored in the vascular network, due to vessels compliance;*

- $\mathcal{K}_j^k(t) = \dfrac{1}{2} L_j^k \left(Q_j^k(t)\right)^2$ *is the kinetic energy related to the blood flow inertia.*

We provide now a deeper explanation of the definition of $\mathcal{E}_\mathrm{j}^\mathrm{k}$. Let us consider, as an example, $\mathcal{E}_\mathrm{AR}^\mathrm{SYS}$. First, we notice that $Q_\mathrm{AV}(t) - Q_\mathrm{AR}^\mathrm{SYS}(t)$ is the net blood flux passing through the arterial systemic network. Hence, by denoting with $V_\mathrm{AR}^\mathrm{SYS}(t)$ the blood volume stored in the arterial systemic network, we have:

$$\frac{dV_\mathrm{AR}^\mathrm{SYS}(t)}{dt} = Q_\mathrm{AV}(t) - Q_\mathrm{AR}^\mathrm{SYS}(t).$$

By comparing the latter equation to the first equation of (7c) we get:

$$p_\mathrm{AR}^\mathrm{SYS}(t) = \frac{1}{C_\mathrm{AR}^\mathrm{SYS}} \left(V_\mathrm{AR}^\mathrm{SYS}(t) - V_\mathrm{0,AR}^\mathrm{SYS}\right),$$

where $V_\mathrm{0,AR}^\mathrm{SYS}$ is the blood volume stored within the arterial systemic network when the pressure is null. In conclusion, we have $\mathcal{E}_\mathrm{AR}^\mathrm{SYS}(t) = (2\,C_\mathrm{AR}^\mathrm{SYS})^{-1} \left(V_\mathrm{AR}^\mathrm{SYS}(t) - V_\mathrm{0,AR}^\mathrm{SYS}\right)^2$, where $1/C_\mathrm{AR}^\mathrm{SYS}$ is the arterial systemic network elastance (the inverse of the compliance), coherently with the definition of $\mathcal{E}_\mathrm{LA}$.

**Definition 2.** *We define the power generated by active contraction, due to ATP consumption occurring at the cellular level, as*

$$\Pi^\mathrm{act}(t) = \Pi_\mathrm{LA}^\mathrm{act}(t) + \Pi_\mathrm{LV}^\mathrm{act}(t) + \Pi_\mathrm{RA}^\mathrm{act}(t) + \Pi_\mathrm{RV}^\mathrm{act}(t),$$

*where* $\Pi_i^\mathrm{act}(t) = -E_i^\mathrm{act}(t)\,(V_i(t) - V_{0,i})\,\frac{dV_i}{dt}(t)$ *is the power exerted by the active contraction of a cardiac chamber (for* $i \in \{\mathrm{LA}, \mathrm{LV}, \mathrm{RA}, \mathrm{RV}\}$*),*

**Definition 3.** *We define the power dissipated within the 0D circulation model by viscous forces as the blood flows through the valves and the vascular network as*

$$\Pi^\mathrm{diss}(t) = \Pi_\mathrm{MV}(t) + \Pi_\mathrm{AV}(t) + \Pi_\mathrm{TV}(t) + \Pi_\mathrm{PV}(t) + \Pi_\mathrm{AR}^\mathrm{SYS}(t) + \Pi_\mathrm{VEN}^\mathrm{SYS}(t) + \Pi_\mathrm{AR}^\mathrm{PUL}(t) + \Pi_\mathrm{VEN}^\mathrm{PUL}(t).$$

*where:*

- *the power dissipated by the blood flux through the cardiac valves:*

$$\begin{aligned} \Pi_\mathrm{MV}(t) &= -\frac{(p_\mathrm{LA}(t) - p_\mathrm{LV}(t))^2}{R_\mathrm{MV}(p_\mathrm{LA}(t), p_\mathrm{LV}(t))}, & \Pi_\mathrm{AV}(t) &= -\frac{\left(p_\mathrm{LV}(t) - p_\mathrm{AR}^\mathrm{SYS}(t)\right)^2}{R_\mathrm{AV}(p_\mathrm{LV}(t), p_\mathrm{AR}^\mathrm{SYS}(t))}, \\ \Pi_\mathrm{TV}(t) &= -\frac{(p_\mathrm{RA}(t) - p_\mathrm{RV}(t))^2}{R_\mathrm{TV}(p_\mathrm{RA}(t), p_\mathrm{RV}(t))}, & \Pi_\mathrm{PV}(t) &= -\frac{\left(p_\mathrm{RV}(t) - p_\mathrm{AR}^\mathrm{PUL}(t)\right)^2}{R_\mathrm{PV}(p_\mathrm{RV}(t), p_\mathrm{AR}^\mathrm{PUL}(t))}; \end{aligned} \tag{11}$$

- $\Pi_j^k(t) = -R_j^k \left(Q_j^k(t)\right)^2$*, that is the power dissipated by the arterial systemic network (for* $j \in \{\mathrm{AR}, \mathrm{VEN}\}$ *and* $k \in \{\mathrm{SYS}, \mathrm{PUL}\}$*).*

We remark that all the terms in Eq. (11) are nonpositive (i.e. dissipative).

**Definition 4.** *We define the power due to the action of the external pressure* $p_{\mathrm{EX}}$ *on the myocardium as*

$$\Pi^{\mathrm{ex}}(t) = \Pi^{\mathrm{ex}}_{\mathrm{LA}}(t) + \Pi^{\mathrm{ex}}_{\mathrm{LV}}(t) + \Pi^{\mathrm{ex}}_{\mathrm{RA}}(t) + \Pi^{\mathrm{ex}}_{\mathrm{RV}}(t),$$

*where* $\Pi^{\mathrm{ex}}_i(t) = -p_{\mathrm{EX}}(t)\dfrac{dV_i}{dt}(t)$ *is the power exerted by the external pressure* $p_{\mathrm{EX}}(t)$ *acting on a cardiac chamber (for* $i \in \{\mathrm{LA}, \mathrm{LV}, \mathrm{RA}, \mathrm{RV}\}$*).*

We have the following result.

**Proposition 1.** *The solution of Eq.* (9) *for the whole 0D circulation model satisfies the energy balance*

$$\frac{d}{dt}\mathcal{M}(t) = \Pi^{\mathrm{act}}(t) + \Pi^{\mathrm{diss}}(t) + \Pi^{\mathrm{ex}}(t), \tag{12}$$

*whose terms are defined in Defs. 1–4.*

*Proof.* To derive Eq. (12) we consider, as illustrative examples, a representative cardiac chamber, a cardiac valve and a vascular branch. In fact, similar calculations will apply to the other chambers, valves and vascular compartments.
*Energy balance of the cardiac chambers.* Let us consider for now the left atrium (LA). By multiplying Eq. (8b) by $\frac{dV_{\mathrm{LA}}(t)}{dt}$ and thanks to Eq. (7a), we get:

$$p_{\mathrm{LA}}(t)(Q^{\mathrm{PUL}}_{\mathrm{VEN}}(t) - Q_{\mathrm{MV}}(t)) = \frac{d}{dt}\mathcal{E}_{\mathrm{LA}}(t) - \Pi^{\mathrm{act}}_{\mathrm{LA}}(t) - \Pi^{\mathrm{ex}}_{\mathrm{LA}}(t). \tag{13}$$

*Energy balance of the cardiac valves.* From (8e) we obtain

$$(p_{\mathrm{LA}}(t) - p_{\mathrm{LV}}(t))Q_{\mathrm{MV}}(t) = -\Pi_{\mathrm{MV}}(t). \tag{14}$$

Similar considerations hold for the other valves.
*Energy balance of the peripheral blood reservoirs.* By multiplying the first equation of (7c) by $p^{\mathrm{SYS}}_{\mathrm{AR}}(t)$, we get:

$$p^{\mathrm{SYS}}_{\mathrm{AR}}(t)Q_{\mathrm{AV}}(t) - p^{\mathrm{SYS}}_{\mathrm{AR}}(t)Q^{\mathrm{SYS}}_{\mathrm{AR}}(t) = \frac{d}{dt}\mathcal{E}^{\mathrm{SYS}}_{\mathrm{AR}}(t) \tag{15}$$

*Energy balance of the peripheral blood conducting system.* By multiplying (7e) by $R^{\mathrm{SYS}}_{\mathrm{AR}}Q^{\mathrm{SYS}}_{\mathrm{AR}}(t)$, we get:

$$p^{\mathrm{SYS}}_{\mathrm{AR}}(t)Q^{\mathrm{SYS}}_{\mathrm{AR}}(t) - p^{\mathrm{SYS}}_{\mathrm{VEN}}(t)Q^{\mathrm{SYS}}_{\mathrm{AR}}(t) = \frac{d}{dt}\mathcal{K}^{\mathrm{SYS}}_{\mathrm{AR}}(t) - \Pi^{\mathrm{SYS}}_{\mathrm{AR}}(t). \tag{16}$$

*Total balance.* By proceeding as above for the other cardiac chambers, valves and circulation systems, and summing up the resulting equations, we obtain Eq. (12). This completes the proof. □

Each of the four terms of Eq. (12) represents the result of the sum of different contributions, associated with the four chambers, the four valves and the different compartments of the vascular network (systemic and pulmonary, arterial and venous). The total work performed in a time interval $[0, T]$ is obtained by integrating the corresponding power over time, according to the following definition.

**Definition 5.** *Let us consider a time horizon* $T > 0$*. The total work performed by active and dissipative forces in the time interval* $[0, T]$ *are defined as*

$$W^{\mathrm{act}} = \int_0^T \Pi^{\mathrm{act}}(t)\, dt, \qquad W^{\mathrm{diss}} = \int_0^T \Pi^{\mathrm{diss}}(t)\, dt,$$

*respectively.*

When the heart is in a periodic regime, it carries out its function alongside a cyclical path. In this case, the work balance of the following proposition holds.

**Proposition 2.** *Let us suppose that $p_{\mathrm{EX}}(t)$ is constant in time. Then, periodic solutions of Eq.* (9) *(i.e. with $\boldsymbol{c}_1(0) = \boldsymbol{c}_1(T)$) satisfy*

$$W^{\mathrm{act}} + W^{\mathrm{diss}} = 0. \tag{17}$$

*Proof.* We integrate the energy balance of Eq. (12) over a cardiac cycle $[0, T]$. Thanks to the periodicity assumption, the contribution of the mechanical energy term $\mathcal{M}$ is null. Moreover, it is easy to show that the term $\Pi^{\mathrm{ex}}(t)$ is conservative and hence its contribution over $[0, T]$ is also null. □

Therefore, when the heart is in a periodic regime, the work performed by the contraction of the four chambers balances the energy dissipated by the four valves and by the blood flux through the systemic and pulmonary circulations.

## 4.2 Energy balance of the 3D-0D coupled model

**Energy balance of the 3D LV model.** By multiplying the first equation of $(\mathscr{M})$ by $\frac{\partial \mathbf{d}}{\partial t}$ and integrating over $\Omega_0$ we obtain:

$$\int_{\Omega_0} \rho_{\mathrm{s}} \frac{\partial^2 \mathbf{d}}{\partial t^2} \cdot \frac{\partial \mathbf{d}}{\partial t} d\mathbf{x} + \int_{\Omega_0} \mathbf{P}(\mathbf{d}, T_{\mathrm{a}}) \colon \nabla \left( \frac{\partial \mathbf{d}}{\partial t} \right) d\mathbf{x} = \int_{\partial\Omega_0} \mathbf{P}(\mathbf{d}, T_{\mathrm{a}}) \mathbf{N} \cdot \frac{\partial \mathbf{d}}{\partial t} \, d\Gamma_0. \tag{18}$$

By substituting the boundary conditions of $(\mathscr{M})$ into (18), we obtain the following energy balance for the 3D left ventricle model:

$$\frac{d}{dt} \mathcal{K}_{\mathrm{LV,3D}}(t) + \frac{d}{dt} \mathcal{E}_{\mathrm{LV,3D}}(t) + \Pi^{\mathrm{act}}_{\mathrm{LV,3D}}(t) + \Pi^{\mathrm{diss}}_{\mathrm{LV,3D}}(t) + \Pi^{\mathrm{press}}_{\mathrm{LV,3D}}(t) = 0. \tag{19}$$

This relation reveals the mutual balance of:

- the kinetic energy associated with the motion of the LV:

$$\mathcal{K}_{\mathrm{LV,3D}}(t) = \frac{1}{2} \int_{\Omega_0} \rho_{\mathrm{s}} \left| \frac{\partial \mathbf{d}}{\partial t} \right|^2 d\mathbf{x};$$

- the elastic energy internally stored by the LV muscle and by the elastic components of the surrounding tissues:

$$\mathcal{E}_{\mathrm{LV,3D}}(t) = \int_{\Omega_0} \mathcal{W}(\mathbf{F}) d\mathbf{x} + \frac{1}{2} \int_{\Gamma_0^{\mathrm{epi}}} \left[ K_{\perp}^{\mathrm{epi}} \left| \mathbf{d} \cdot \mathbf{N} \right|^2 + K_{\parallel}^{\mathrm{epi}} \left| (\mathbf{I} - \mathbf{N} \otimes \mathbf{N}) \mathbf{d} \right|^2 \right] d\Gamma_0,$$

  where the displacement at the epicardium is split into the normal $|\mathbf{d} \cdot \mathbf{N}|$ and tangent $|(\mathbf{I} - \mathbf{N} \otimes \mathbf{N})\mathbf{d}|$ component;

- the power exerted by the active contraction of the LV:

$$\Pi^{\mathrm{act}}_{\mathrm{LV,3D}}(t) = - \int_{\Omega_0} T_{\mathrm{a}} \frac{\mathbf{F} \mathbf{f}_0 \otimes \mathbf{f}_0}{\sqrt{\mathcal{I}_{4f}}} : \nabla \left( \frac{\partial \mathbf{d}}{\partial t} \right) d\mathbf{x},$$

- the power dissipated by the interaction with the pericardium:

$$\Pi^{\text{diss}}_{\text{LV,3D}}(t) = -\int_{\Gamma_0^{\text{epi}}} \left[ C_{\perp}^{\text{epi}} \left| \frac{\partial \mathbf{d}}{\partial t} \cdot \mathbf{N} \right|^2 + C_{\parallel}^{\text{epi}} \left| (\mathbf{I} - \mathbf{N} \otimes \mathbf{N}) \frac{\partial \mathbf{d}}{\partial t} \right|^2 \right] d\Gamma_0 \leq 0,$$

  $\forall t$, which is a nonpositive term;

- the power exchanged with the blood contained in the LV cavity, by means of the action of pressure $p_{\text{LV}}(t)$ on the endocardium:

$$\Pi^{\text{press}}_{\text{LV,3D}}(t) = p_{\text{LV}}(t) \left[ \int_{\Gamma_0^{\text{endo}}} J\mathbf{F}^{-T}\mathbf{N} \cdot \frac{\partial \mathbf{d}}{\partial t} \, d\Gamma_0 - \int_{\Gamma_0^{\text{base}}} |J\mathbf{F}^{-T}\mathbf{N}| \frac{\partial \mathbf{d}}{\partial t} \, d\Gamma_0 \cdot \mathbf{v}^{\text{base}}(t) \right].$$

  As in [44], the term in square brackets corresponds to the time derivative of the LV volume, that is $\Pi^{\text{press}}_{\text{LV,3D}}(t) = p_{\text{LV}}(t) \frac{d}{dt} V^{\text{3D}}_{\text{LV}}(\mathbf{d}(t))$. Then, in virtue of the coupling condition $(\mathscr{V})$ and by Eq. (7a), we obtain $\Pi^{\text{press}}_{\text{LV,3D}}(t) = p_{\text{LV}}(t)(Q_{\text{MV}}(t) - Q_{\text{AV}}(t))$.

**Total energy balance.** By setting $p_{\text{EX}}(t) \equiv 0$, i.e. by neglecting the effect of the pressure exerted by the surronding organs, and by replacing the energy balance of the 0D left ventricle model with Eq. (19), we obtain again Eq. (12), where the total mechanical energy is now

$$\begin{aligned} \mathcal{M}(t) = {} & \mathcal{E}_{\text{LA}}(t) + \mathcal{E}_{\text{LV,3D}}(t) + \mathcal{E}_{\text{RA}}(t) + \mathcal{E}_{\text{RV}}(t) + \mathcal{E}^{\text{SYS}}_{\text{AR}}(t) + \mathcal{E}^{\text{SYS}}_{\text{VEN}}(t) + \mathcal{E}^{\text{PUL}}_{\text{AR}}(t) + \mathcal{E}^{\text{PUL}}_{\text{VEN}}(t) \\ & + \mathcal{K}^{\text{SYS}}_{\text{AR}}(t) + \mathcal{K}^{\text{SYS}}_{\text{VEN}}(t) + \mathcal{K}^{\text{PUL}}_{\text{AR}}(t) + \mathcal{K}^{\text{PUL}}_{\text{VEN}}(t) + \mathcal{K}_{\text{LV,3D}}(t) \end{aligned} \tag{20}$$

and the power of active contraction and the total dissipated power ($\Pi^{\text{diss}}(t) \geq 0$) are

$$\begin{aligned} \Pi^{\text{act}}(t) &= \Pi^{\text{act}}_{\text{LA}}(t) + \Pi^{\text{act}}_{\text{LV,3D}}(t) + \Pi^{\text{act}}_{\text{RA}}(t) + \Pi^{\text{act}}_{\text{RV}}(t); \\ \Pi^{\text{diss}}(t) &= \Pi_{\text{MV}}(t) + \Pi_{\text{AV}}(t) + \Pi_{\text{TV}}(t) + \Pi_{\text{PV}}(t) \\ &\quad + \Pi^{\text{SYS}}_{\text{AR}}(t) + \Pi^{\text{SYS}}_{\text{VEN}}(t) + \Pi^{\text{PUL}}_{\text{AR}}(t) + \Pi^{\text{PUL}}_{\text{VEN}}(t) + \Pi^{\text{diss}}_{\text{LV,3D}}(t), \end{aligned}$$

respectively; here $\Pi^{\text{ex}}(t) \equiv 0$. We remark that Prop. 2 applies also to this case. Finally, we conclude that the 3D-0D coupling is compliant with the princicple of conservation of mechanical energy. We remark that this result is achieved thanks to the energy-consistent boundary conditions imposed at the left ventricle base; see $(\mathscr{M})$. As a matter of fact, if other boundary conditions – such as homogeneous Neumann conditions – are imposed at the base instead, the relationship $\Pi^{\text{press}}_{\text{LV,3D}}(t) = p_{\text{LV}}(t) \frac{d}{dt} V^{\text{3D}}_{\text{LV}}(\mathbf{d}(t))$ may not hold and the balance of Eq. (12) is not satisfied.

We notice that, compared to the fully 0D case, the 3D electromechanical model shows two additional terms, namely $\mathcal{K}_{\text{LV,3D}}(t)$ and $\Pi^{\text{diss}}_{\text{LV,3D}}$, respectively accounting for the kinetic energy of the LV and for the dissipation associated with the interaction of $\Gamma_0^{\text{epi}}$ with surrounding tissues. Indeed, both features are not included in the 0D circulation model, in which cardiac chambers are modeled quasistatically.

## 4.3 Quantitative analysis of cardiac energetics

In Fig. 4 we report all the energy and power terms over a characteristic steady-state cardiac cycle obtained with a numerical simulation of the 0D circulation model. For this simulation we employ the parameters reported in Part II of this paper [45], with $E^{\text{act,max}}_{\text{LV}} = 2.75$ mmHg mL$^{-1}$ and $E^{\text{pass}}_{\text{LV}} = 0.08$ mmHg mL$^{-1}$. Fig. 4 (top-left) displays the time evolution of the terms of Eq. (12). We notice that, while the energy input ($\Pi^{\text{act}}$) occurs in a short time interval of nearly 100 ms (during systole), energy dissipation ($\Pi^{\text{diss}}$) takes place throughout the entire duration of the heartbeat. As a matter of fact, mechanical energy $\mathcal{M}$ plays a dominant role. Moreover, it is initially accumulated and then it is gradually dissipated as the blood flows through systemic and pulmonary circulations.

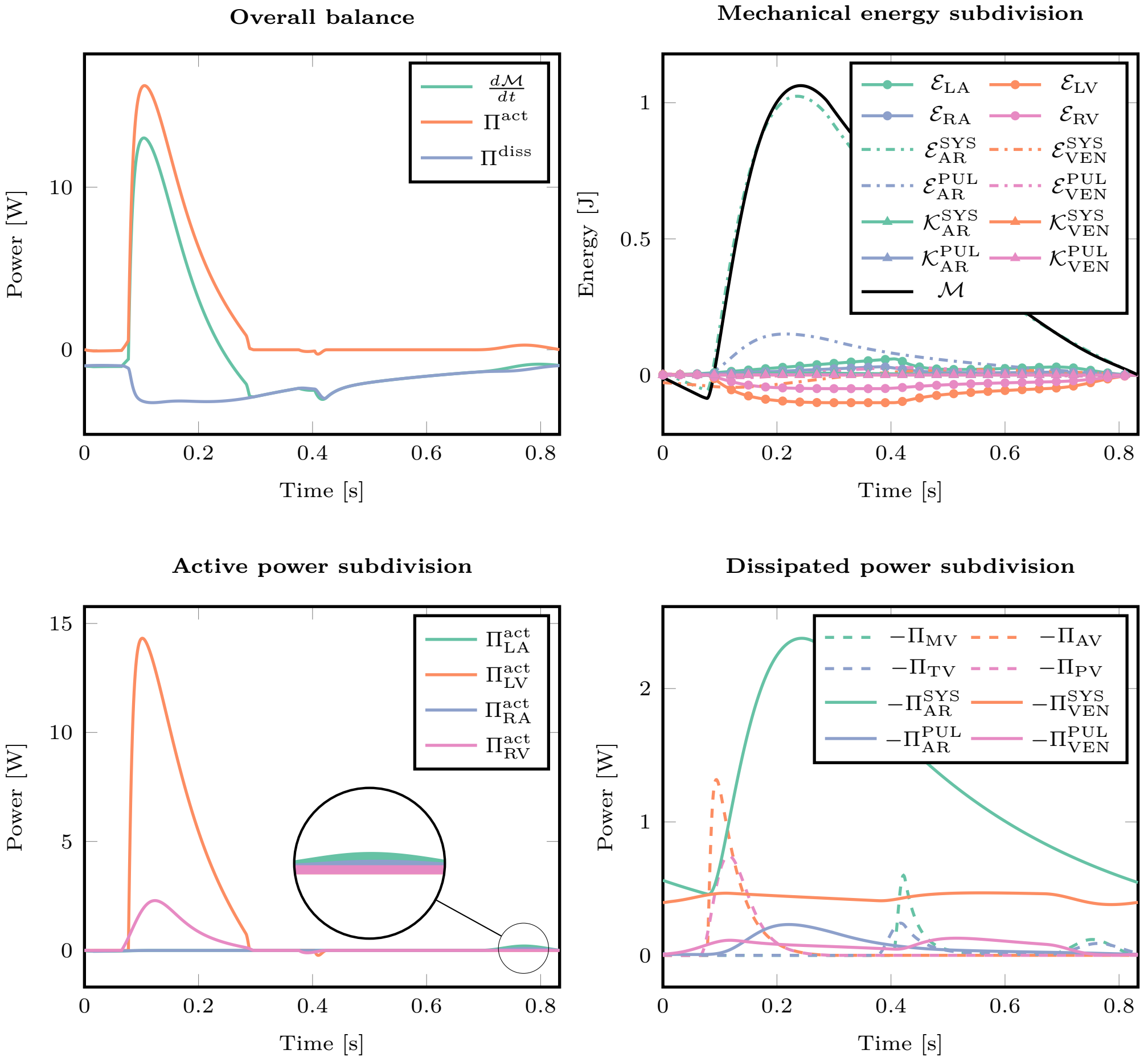

Figure 4: Time evolution of both power and energy terms ($\mathcal{M}$, $\mathcal{E}$, $\mathcal{K}$, $\Pi^{\text{act}}$, $\Pi^{\text{diss}}$) of the 0D circulation model. A single heartbeat in a periodic regime is considered.

Fig. 4 (top-right, bottom-left, bottom-right) illustrate the details of the three terms $\mathcal{M}$, $\Pi^{\text{act}}$ and $\Pi^{\text{diss}}$, showing how they are divided into the various subterms during the different phases of the heartbeat. Specifically, we notice that the chamber that contributes the most to the work generation is the LV, followed by the RV, while the atria only contribute – albeit to a small extent – around $t = 0.8$ s, during the atrial systole. The large part of mechanical energy and of dissipated power are associated with the systemic arterial circulation, as this branch of the circulatory network is located downstream the LV, the chamber carrying out most of the mechanical work. We remark that a non negligible dissipation of energy also takes place across the open valves, due to the high-speed blood flow across the valvular orifices.

Our model allows to estimate the daily production of mechanical work of the heart. This is obtained by multiplying the number of seconds in a day times the average generated power, given by:

$$\overline{\Pi}^{\mathrm{act}} = \frac{1}{T} \int_0^T \Pi^{\mathrm{act}}(t)\, dt.$$

Applying this formula to the results of the simulation considered in Fig. 4, we obtain a daily work production of 182.5 kJ, of which 155.9 kJ attributable to the left ventricle, 24.8 kJ to the right ventricle and only 1.8 kJ to the atria.

In the daily clinical practice, the work generated by the myocardium is instead estimated through simple relationships [25, 35]. In this regard, our model offers a tool to estimate the validity of these approaches. A commonly used formula [25] is

$$\overline{\Pi}^{\mathrm{act}} \simeq \overline{p_{\mathrm{AR}}^{\mathrm{SYS}}} \frac{SV}{T_{\mathrm{beat}}}, \tag{21}$$

where $\overline{p_{\mathrm{AR}}^{\mathrm{SYS}}}$ denotes the average systemic arterial pressure (corresponding to the wrist average blood pressure), $SV$ is the stroke volume (i.e. the difference between maximum and minimum $V_{\mathrm{LV}}$) and $T_{\mathrm{beat}}$ is the heartbeat period. By applying (21) to the results of the above simulation, we obtain a daily work generation of 152.4 kJ. Hence (21) underestimates the mechanical work by 16%. As a matter of fact, (21) only refers to the work done by the left ventricle (which is, instead, approximated up to an error of only 2%). The large part of the error is thus attributable to the work performed by the right ventricle, which is not accounted for in (21).

In common clinical practice, then, $\overline{p_{\mathrm{AR}}^{\mathrm{SYS}}}$ is not directly measured, but it is estimated as $1/3\, p_{\max} + 2/3\, p_{\min}$, where $p_{\max}$ and $p_{\min}$ are the maximum (systolic) and minimum (diastolic) arterial pressures (see e.g. [35]). With this further approximation we obtain an estimated daily work of 138.8 kJ: this underestimates the left ventricle work by 11% and the total work of the myocardium by 24%.

# 5 Conclusions

We presented a mathematical model of cardiac electromechanics, where the different physical phenomena therein involved are described by means of biophysically detailed core models. To provide a realistic and physically meaningful relationship between the blood flux through the left ventricle and the pressure acting against the endocardial surface, we coupled the 3D electromechanical model of the left ventricle with a closed-loop 0D model of the whole circulation. We proved that our 3D-0D closed-loop model of the whole circulation is compliant with the principles of conservation of mechanical energy. Indeed, the power exerted by the cavity pressure in the 3D electromechanical model balances the power exchanged with the 0D circulation model at the coupling interface. This is to be ascribed to the energy-consistent boundary conditions, that we proposed in [44], adopted in the 3D mechanical problem.

In our active stress model we replaced a biophysically detailed, but computationally demanding, subcellular model of active force generation by a surrogate model, based on an ANN [40]. This model, built by means of Machine Learning from a collection of pre-computed simulations, allows to accurately reproduce the results of the high-fidelity model by reducing by a factor of 1000 the number of internal variables [43, 44]. In this way we obtain in a very favorable trade-off between the biophysical accuracy of the results and the computational cost of our numerical simulations.

We have formulated an inverse problem aimed at recovering the stress-free configuration from a stressed geometry, which, in practical applications, corresponds to the geometry acquired from medical imaging. After this reference configuration is obtained, we can inflate the left ventricle by applying the end-diastolic pressure at the endocardium, thus obtaining the initial displacement needed to start the numerical simulation.

We analyzed the mechanical work associated with the different compartments of our circulation model and we proved that a balance of mechanical energy is satisfied. This balance holds both when we consider the 0D circulation model alone and when we consider the 3D-0D coupled model. We showed that the circulation model considered in this paper can be exploited to provide quantitative insight into the heart energy distribution. In particular we employed this model to validate the reliability of relationships used in the daily clinical practice to estimate the mechanical work performed by the heart. We highlighted that these relationship can be accurate when used to assess the left ventricle function, but less accurate when the mechanical work of the whole myocardium is addressed.

## Acknowledgements

This project has received funding from the European Research Council (ERC) under the European Union's Horizon 2020 research and innovation programme (grant agreement No 740132, iHEART - An Integrated Heart Model for the simulation of the cardiac function, P.I. Prof. A. Quarteroni). We acknowledge the CINECA award under the class C ISCRA project HP10CWQ2GS, for the availability of high performance computing resources and support.